\documentclass[11pt]{article} 

\usepackage{amsmath}
\usepackage{amssymb}
\usepackage{color}
\usepackage[T1]{fontenc}
\usepackage[latin1]{inputenc}
\def\r{\rightarrow}

\usepackage{enumerate}

\newcommand{\fdem}{\hspace*{\fill}~$\Box$\par\endtrivlist\unskip}

\newcommand{\E}{\mathbb{E}}     
     
\renewcommand{\L}{\mathbb{L}}

\newcommand{\N}{\mathbb{N}}     

\newcommand{\R}{\mathbb{R}}     
     
\newcommand{\C}{\mathbb{C}}

\newcommand{\X}{\mathbb{X}}

\renewcommand{\Im}{\mbox{\rm Im }}

\renewcommand{\dim}{\mathop{\rm dim}}
\renewcommand{\ker}{\mathop{\rm Ker}}

\renewcommand{\r}{\mathop{\rightarrow}}

\newcommand{\cB}{\mbox{$\cal B$}}

\newcommand{\cD}{\mbox{$\cal D$}}

\newcommand{\cL}{\mbox{$\cal L$}}

\newcommand{\cX}{\mbox{$\cal X$}}

\newtheorem{theo}{Theorem}
\newtheorem{pro}{Proposition}
\newenvironment{proof}[1]{\textit{Proof#1.\,}}{\fdem}
\newtheorem{lem}{Lemma}
\newtheorem{rem}{Remark}
\newtheorem{cor}{Corollary}

\newtheorem{ex}{Example}

\title{Regular perturbation of $V$-geometrically ergodic Markov chains}
\author{Déborah FERR\'E, Loïc HERV\'E, James LEDOUX \footnote{Université européenne de Bretagne, France ; IRMAR UMR CNRS 6625;  Institut National des Sciences Appliquées de Rennes; Deborah.Ferre,Loic.Herve,James.Ledoux@insa-rennes.fr}
}

\begin{document}
\maketitle

\begin{abstract}
In this paper, new conditions for the stability of $V$-geometrically ergodic Markov chains are introduced. 
The results are based on an extension of the standard perturbation theory formulated by Keller and Liverani. The continuity and higher regularity properties are investigated. As an illustration, an asymptotic expansion of the invariant probability measure for an autoregressive model with i.i.d.~noises (with a non-standard probability density function) is obtained.
\end{abstract}

\begin{center}


\end{center}
\section{Introduction and statements}

Let  $\{P_\varepsilon\}_{|\varepsilon|<\varepsilon_0}$ be a family of transition kernels on a measurable space $(\X,\cX)$, where $\varepsilon$ reads as a small perturbation parameter. Throughout the paper $V : \X\r [1,+\infty)$ is a fixed function. The (unperturbed) kernel $P_0$ is assumed to satisfy the classical $V$-geometrical  ergodicity property, called (\ref{VG}) in our paper,
that is: $P_0$ admits a unique invariant probability measure $\pi_0$, we have $\pi_0(V)<\infty$, and there exists some constants $c\in(0,+\infty)$ and $ \kappa_1\in(0,1)$ such that: 
\begin{equation} \label{VG}
\forall x \in\X, \quad \sup_{|f|\leq V} \big|\, \E[f(X_n)\, |\, X_0=x] - \pi_0(f)\, \big| \leq c\,  \kappa_1^n\, V(x).  \tag{\bf VG} 
\end{equation} 
or  equivalently 
\begin{equation*} 
\forall x \in\X, \quad \sup_{|f|\leq V} \big|P_0^n f(x)-\pi_0(f)\big| \leq c\,  \kappa_1^n\, V(x).
\end{equation*} 
This means that $P_0$ has a spectral gap on the weighted supremum normed space $\cB_V$ composed of the measurable functions $f : \X\r\C$ such that $\|f\|_V:=\sup_{x\in\X} V(x)^{-1}|f(x)| < \infty$. We are interested in the two following questions. For $|\varepsilon|$ small enough,\textit{
\begin{enumerate}[(I)] 
	\item  \label{I} does $P_\varepsilon$ admit an invariant probability measure, say $\pi_\varepsilon$, and is $P_\varepsilon$ $V$-geometrically ergodic?
	\item  \label{II} Do we have any control on $\pi_0-\pi_\varepsilon$?
\end{enumerate}}

Under some classical aperiodicity and irreducibility conditions, the  property (\ref{VG}) holds true if and only if $P_0$ satisfies the so-called drift condition based on the notion of small set (see \cite{mey} for the definition of the drift condition which is not used here). Consequently a natural and efficient way to study {\it(\ref{I})} is to prove that the perturbed Markov kernel $P_\varepsilon$ also satisfies the drift condition (e.g.~see \cite{rob-ros-sch,bre-rob-ros}). However, to the best of our knowledge, the theory of geometrical ergodic Markov chains does not provide any general answer for question {\it(\ref{II})}, except in terms of weak convergence in specific cases (see \cite{rob-ros-sch}). On the other hand the standard perturbation theory, which is a natural way to investigate {\it(\ref{II})}, leads to assume the following continuity condition: 
\begin{equation} \label{Cont_11} 
\|P_{\varepsilon} - P_0 \|_{{\cal B}_V} := \sup _{\|f\|_V\le 1} \| P_\varepsilon f- P_0f\|_V \r0\ \text{ when  } \varepsilon\r 0
\end{equation} 
that is the operator norm of $P_{\varepsilon} - P_0$ on $\cB_V$  goes to 0 as $\varepsilon\r 0$. 
In a series of papers, Kartashov has introduced the concept of ``strongly stable Markov chain'' for a Markov chain with a transition kernel $P_0$ such that, in some neighborhood of $P_0$ with respect to $\|\cdot\|_{{\cal B}_V}$, $P_\varepsilon$ has a unique invariant probability measure $\pi_{\varepsilon}$ with the property $\sup_{\|f\|_V\le 1} | \pi_{\varepsilon}(f) -\pi_0(f) |\r 0$ as $\| P_{\varepsilon} - P_0 \|_{{\cal B}_V} \r 0$ uniformly in this neighborhood. Strong stability is shown to be equivalent to the convergence  
\begin{equation*} 
\lim_{n \r \infty}\sup _{\|f\|_V\le 1} \big\| \frac{1}{n}\sum_{k=1}^n P_0^k f- \pi_0(f)1_{\X}\big\|_V=0.
\end{equation*} 
Moreover, if $P_0$ is $V$-geometrically ergodic, then given $\rho\in(0,1)$, one can consider $N\in\N^*$ such that $c\kappa_1^N\leq\rho$, and next if $P_\varepsilon$ is such that $\Delta_N := \| P_0^N - P_\varepsilon^N \|_{{\cal B}_V} < 1-\rho$, then $P_\varepsilon$ is $V$-geometrically ergodic and $\sup_{\|f\|_V\le 1} | \pi(f) - \pi_{\varepsilon}(f)| = O(\Delta_N/(1-\rho-\Delta_N))$. We refer to \cite{kar96} for an overview of results in this direction and to \cite{RabAis10} for a related discussion for discrete state spaces. However, as discussed in \cite[p.~1126]{shar-stua} and in Example~\ref{AR_weak_pasStrong} below, the continuity condition (\ref{Cont_11}) may be 
restrictive. 

Similar questions arise in the context of dynamical systems. To overcome the previous difficulty, Keller introduced the more general assumption 
$$ \lim_{\varepsilon\r 0}\sup_{\|f\|_0\leq1}\|P_\varepsilon f- P_0 f\|_1 = 0  $$
 involving two norms $\|\cdot\|_0$ and $\|\cdot\|_1$ (instead of a single one) on the space on which $P_0$ has a spectral gap \cite{kel82}. This approach has been highly enhanced by the Keller-Liverani perturbation theorem \cite{kelliv99,liv04}, which has proved to be very powerful in studying the behaviour of the  Sinai-Ruelle-Bowen measures of certain perturbed dynamical systems (e.g.~see \cite[Th 2.10]{bal} and \cite[Th.~2.8]{gouliv}). 
 
The goal of this paper is to show that the Keller-Liverani theorem also provides an interesting way to investigate both the questions {\it(\ref{I})} {\it(\ref{II})} in the context of geometrical ergodic Markov chains. In this markovian context, the closest work to ours is \cite{shar-stua} where Keller's approach is used. The results of \cite{shar-stua} are improved here thanks to the Keller-Liverani perturbation theorem. Furthermore in this paper, higher regularity properties than continuity are investigated in question~{\it(\ref{II})}. Mention that the results of  \cite{kelliv99,liv04} have been  already used in \cite{deb-2011} to study some stability properties of parametric autoregressive models  (for  different purposes from those of Proposition~\ref{auto-reg} below). 

\paragraph{Notations.} For $\beta\in[0,1]$, we denote by $(\cB_{\beta},\|\cdot\|_{\beta})$ the Banach space composed of the  measurable functions $f : \X\r\C$ such that $\|f\|_{\beta}  := \sup_{x\in\X} V(x)^{-\beta}|f(x)| < \infty$. Note that $\cB_0$ corresponds to the space of bounded measurable functions on $\X$, with $\|f\|_0 = \sup_{x\in\X}|f(x)|$, and that $\cB_1 =\cB_V$. We denote by $(\cL(\cB_{\beta},\cB_{\beta'}),\|\cdot\|_{\beta,\beta'})$ the space of all the bounded linear maps from $\cB_{\beta}$ to $\cB_{\beta'}$, equipped with its usual norm: $\|T\|_{\beta,\beta'} = \sup\big\{\|Tf\|_{\beta'},\, f\in\cB_\beta,\, \|f\|_\beta \leq 1\big\}$. 
We write $\cL(\cB_\beta)$ for $\cL(\cB_\beta,\cB_\beta)$ and $\|T\|_\beta$ for $\|T\|_{\beta,\beta}$ which is a slight abuse of notation. Let $(\cB_{\beta}',\|\cdot\|_{{\cal B}_\beta'})$ denote the dual space of $\cB_\beta$. 
If $T\in\cL(\cB_{\beta})$, then $T^*$ stands for the adjoint operator of $T$. By definition we have $T^*\in\cL(\cB_{\beta}')$ with the corresponding operator norm $\| T^*\|_{{\cal B}'_{\beta}} = \| T\|_{\beta}$. Note that $T^*$ also defines an element of $\cL(\cB_{\beta}',\cB_0')$ with corresponding operator norm $\| T^*\|_{{\cal B}'_{\beta},{\cal B}'_0 } \leq \| T^*\|_{{\cal B}'_{\beta}}$ from the continuous inclusion ${\cal B}_0\subset{\cal B}_{\beta}$. 

 Each perturbed Markov kernel $P_\varepsilon$ is assumed to continuously act on $\cB_1$. The unperturbed kernel $P_0$ is assumed to satisfy (\ref{VG}), namely: $P_0$ admits a unique invariant distribution $\pi_0$ on $(\X,\cX)$, $\pi_0(V) < \infty$, and  
\begin{equation} \label{V-geo-cond} 
\exists  \kappa_1\in(0,1),\quad \|P_0^n - \pi_0(\cdot)1_\X\|_{1} = O( \kappa_1^n). \tag{$V_1$}
\end{equation}
We also assume that there exist $N\in\N^*$, $L\in(0,+\infty)$ and $\delta\in(0,1)$ such that 
\begin{equation} \label{drift-gene-cond}
\forall \varepsilon\in(-\varepsilon_0,\varepsilon_0),\quad P_\varepsilon^N V \leq \delta^N V + L\, 1_{\X}. \tag{D}
\end{equation}
Under some classical aperiodicity and irreducibility assumptions, Property~(\ref{V-geo-cond}) is equivalent to the drift condition (see \cite[Chapter 4, Th. 16.0.1]{mey}).  Condition~(\ref{drift-gene-cond}) on the family $\{P_\varepsilon\}_{|\varepsilon|<\varepsilon_0}$ is weaker than the simultaneous geometrical ergodicity condition introduced in \cite{rob-ros-sch} since (\ref{drift-gene-cond}) involves no small set.\footnote{In case $\varepsilon=0$ (i.e.~for a single transition kernel $P_0$), the connection between the ``best'' constants $\kappa_1$ in (\ref{V-geo-cond}) and $\delta$ in~(\ref{drift-gene-cond}) is discussed in details in \cite{djl-qc}. See also Remark~\ref{Rem_papier_qc}.} 
\begin{theo} \label{pro-bv}
Conditions {\em (\ref{V-geo-cond})}, {\em (\ref{drift-gene-cond})} and 
\begin{equation} \label{Cont_01}
\lim_{\varepsilon\r0} \|P_\varepsilon  - P_0\|_{0,1} = 0
\end{equation}
 are assumed to hold. 
Then, setting $\widehat\kappa := \max(\kappa_1,\delta)$, the following statements are fulfilled: 
\begin{enumerate}
	\item for each $\kappa\in(\widehat\kappa,1)$, there exists $\varepsilon_1\in(0,\varepsilon_0]$ such  that, for all $\varepsilon\in(-\varepsilon_1,\varepsilon_1)$, $P_\varepsilon$ has a unique invariant probability measure $\pi_\varepsilon$, with $\pi_\varepsilon(V)<\infty$, such that 
\begin{equation} \label{unif-bound-rate}
\sup_{|\varepsilon|<\varepsilon_1}\|P_\varepsilon^n - \pi_\varepsilon(\cdot)1_\X\|_{1} = O(\kappa^n).
\end{equation}
	\item We have 
\begin{equation} \label{pro-bv-cont}
	\lim_{\varepsilon\r0} \sup_{\|f\|_0\leq1} \big|\pi_\varepsilon(f) - \pi_0(f)\big| = 0.
\end{equation}	
\end{enumerate}
\end{theo}
Note that, if each probability measure $\pi_\varepsilon$ admits a density, say $p_\varepsilon$, with respect to a fixed positive measure $\psi$ on $(\X,\cX)$, then we have $\sup_{\|f\|_0\leq1} |\pi_\varepsilon(f) - \pi_0(f)| = \int_\R |p_{\varepsilon}(x) - p_{0}(x)|\, d\psi(x)$. 

Conclusion~(\ref{unif-bound-rate}) means that $P_\varepsilon$ is $V$-geometrically ergodic in a uniform way with respect to the perturbation parameter $\varepsilon$. Conclusion~(\ref{pro-bv-cont}) means that the total variation norm of $\pi_\varepsilon - \pi_0$ goes to 0 when $\varepsilon\r0$. Consequently Theorem~\ref{pro-bv} provides the same theoretical conclusions that in \cite[Chap.~3]{kar96}, but under the continuity condition (\ref{Cont_01}) which is weaker than that in (\ref{Cont_11}). Indeed we have: $\|P_\varepsilon  - P_0\|_{0,1} \leq \|P_\varepsilon  - P_0\|_{1}$. There are several examples (see for instance \cite{shar-stua}) showing that the continuity condition (\ref{Cont_11}) of \cite{kar96} may fail, while (\ref{Cont_01}) (or (\ref{Cont_beta_1}) below) holds true. The case of an AR process is investigated in Example~\ref{AR_weak_pasStrong} below. Note that the uniform bound~(\ref{unif-bound-rate}) can be obtained from \cite{mey-bound} when the kernels $P_\varepsilon$ are assumed to satisfy the drift condition with  constants and small set which do not depend on $\varepsilon$. To derive (\ref{unif-bound-rate}) and (\ref{pro-bv-cont}) from the weaker assumptions (\ref{V-geo-cond}) and (\ref{drift-gene-cond}), some continuity assumption must be assumed on the map $\varepsilon\mapsto P_\varepsilon$.  

Typically Condition~(\ref{V-geo-cond}) is the assumption~AI in \cite{shar-stua} (see the remark following their assumptions). Condition~(\ref{Cont_01}) is weaker than the assumption AII in \cite{shar-stua} involving the sequence $(P_\varepsilon^n)_{n\ge 1}$ (here Assumption~AII is only required for $n=1$). 
Mention that, given any $r\in(\widehat\kappa,1)$ and setting $\eta := 1 - \ln r/\ln \delta\, $ ($\eta\in(0,1)$), the proof of Theorem~\ref{pro-bv} and further results in \cite{kelliv99} ensure the following: 
\begin{equation} \label{raf-KL}
\exists D_r\in(0,+\infty), \quad \sup_{\|f\|_0\leq1} \big|\pi_\varepsilon(f) - \pi_0(f)\big|  \leq D_r\, (\|P_\varepsilon  - P_0\|_{0,1})^\eta.
\end{equation}
This provides an alternative statement to \cite[Th.~3.1]{shar-stua}, which states an inequality of type (\ref{raf-KL}) in assuming the existence of $\pi_\varepsilon$. The $V$-geometrical ergodicity of the perturbed kernel $P_\varepsilon$ was an open question in \cite{shar-stua}. Actually, under their assumptions, Condition~(\ref{drift-gene-cond}) is a quite natural hypothesis for $P_\varepsilon$ to inherit the $V$-geometrical ergodicity of the unperturbed kernel $P_0$ (note that Condition~(\ref{drift-gene-cond}) is not so far from the drift condition). In the proof of Theorem~\ref{pro-bv}, Condition~(\ref{drift-gene-cond}) is viewed as a Doeblin-Fortet inequality on the dual space of $\cB_1$ in order to use the Keller-Liverani theorem. 
\begin{ex}[AR process] \label{AR_weak_pasStrong}

Assume that $\X:=\R$ and $(X_n)_{n\in\N}$ is the autoregressive model defined by 
\begin{equation} \label{def-AR}
n\in\N^*, \quad X_n = \alpha X_{n-1} + \vartheta_n\, 
\end{equation}
where $X_0$ is a real-valued random variable, $\alpha\in(-1,1)$, and $(\vartheta_n)_{n\ge 1}$ is a sequence of  i.i.d.~real-valued random variables, independent of $X_0$.  Assume that $\vartheta_1$ has a Lebesgue probability density function on $\X$, say $\nu(\cdot)$, 
and admits a first moment, $\int |x|\nu(x) dx <\infty$. 
We know that $(X_n)_{n\in\N}$ is a Markov chain with transition kernel 
\begin{equation} \label{noyau_AR}
P_{\alpha}(x,A)=\int_{\R} 1_A(\alpha x + y) \nu(y)dy = \int_{\R} 1_A( y) \nu(y-\alpha x)dy .
\end{equation}
Set $V(x) := 1+|x|$, $x\in\R$.
It is known (e.g.~\cite[Sec. 8]{Wu04} or \cite[Sec. 5.5]{djl-qc}) that, for each $\alpha\in(-1,1)$, $(X_n)_{n\in\N}$ is $V$-geometrically ergodic with an invariant distribution $\pi_\alpha$. Next, given any $a_0\in(0,1)$, it can be easily checked that the family $\{P_\alpha,\, \alpha\in(-a_0,a_0)\}$ satisfies Condition~{\em (\ref{drift-gene-cond})} with $N=1$ and with any $\delta\in(a_0,1)$. Finally we prove below that, for every $\alpha_0\in(-1,1)$, the kernels $P_{\alpha_0+\varepsilon}$ satisfy the weak continuity condition~(\ref{Cont_01}). The previous facts ensure that Theorem~\ref{pro-bv} applies to the family $(P_{\alpha_0+\epsilon})_\epsilon$, so that the total variation norm of $\pi_\alpha - \pi_{\alpha_0}$ goes to 0 when $\alpha\r\alpha_0$. Note that no specific assumptions are required on the density $\nu$, excepted that it admits a first moment. To the best of our knowledge this last result is new. Anyway, it cannot be deduced from a known method as in \cite{kar96} since Condition~(\ref{Cont_11}) does not hold, as proved afterwards. 
 
First we show that, for any $\alpha_0\in(-1,1)$, we have $\lim_{\alpha\r\alpha_0} \|P_{\alpha}  - P_{\alpha_0}\|_{0,1} = 0$, that is  
	\[ \forall \varepsilon >0,\  \exists \eta >0, \quad |\alpha -\alpha_0|<\eta \Longrightarrow \| P_{\alpha}-P_{\alpha_0}\|_{0,1} = \sup_{\|f\|_0\leq1} \sup_{x\in\R}\frac{|P_\alpha f(x)- P_{\alpha_0} f(x) |}{V(x)} < \varepsilon .
\]
From (\ref{noyau_AR}), we have $P_{\alpha} f(x) = \int_{\R} f( y) \nu(y-\alpha x)dy$ so that 
	\[ \| P_{\alpha}-P_{\alpha_0}\|_{0,1} \le  \sup_{x\in\R}\frac{\|\nu(\cdot -\alpha x) - \nu(\cdot -\alpha_0 x)  \|_{\L^ 1(\R)}}{V(x)}. 
\]
Let $\varepsilon >0$. Since $V(x) \r +\infty$ when $x\r \infty$, we can choose a constant $A>0$ such that $4/V(x) < \varepsilon$ for any $|x|> A$. Therefore, using that $\nu$ is a probability density function, we obtain   
\[  \sup_{|x| >A}\frac{\|\nu(\cdot -\alpha x) - \nu(\cdot -\alpha_0 x)  \|_{\L^ 1(\R)}}{V(x)} \le 
\sup_{|x| >A}\frac{2}{V(x)} <\varepsilon.\]
 Next, the map from $\R$ into $\L^1(\R)$ defined by $t\mapsto \nu(\cdot - t)$ is continuous, so that we can pick $\eta' >0$ such that 
	\[\forall |x|\le A, \quad |\alpha x-\alpha_0x| \le |\alpha-\alpha_0| A < \eta' \Longrightarrow \| \nu(\cdot - \alpha x)- \nu(\cdot - \alpha_0 x)\|_{\L^ 1(\R)} <\varepsilon.
\]
Therefore, given that $V\ge 1$, there exists $\eta:=\eta'/A$ such that 
 \[ |\alpha-\alpha_0| < \eta \Longrightarrow \sup_{|x|\le A}\frac{\|\nu(\cdot -\alpha x) - \nu(\cdot -\alpha_0 x)  \|_{\L^ 1(\R)}}{V(x)}  
 \le  \sup_{|x|\le A}\|\nu(\cdot -\alpha x) - \nu(\cdot -\alpha_0 x)  \|_{\L^ 1(\R)}<\varepsilon.\]

Second, let us check that, whatever the density $\nu(\cdot)$ is, the strong continuity condition~(\ref{Cont_11}) is never fulfilled. Let us consider a positive real number $a$ such that 
$$\int_{-a}^a\nu(y)dy \neq \int_{-2a}^{-a}\nu(y)dy + \int_a^{2a}\nu(y)dy.$$
Such a real number exists for any density $\nu(\cdot)$ (if the previous terms coincide for every $a>0$, then we obtain a contradiction when $a\r+\infty$). Let $\alpha_0\in(0,1)$ be fixed. For each $\alpha\in(\alpha_0,1)$, define 
$$x_\alpha := \frac{a}{\alpha-\alpha_0}>0.$$
Note that $-a+\alpha_0\, x_\alpha < -a+\alpha\, x_\alpha < a + \alpha_0\, x_\alpha < a + \alpha\, x_\alpha$. Next let us introduce the following function 
$$\forall y\in\R,\quad f_\alpha(y) := y\, 1_{ [a+\alpha_0 x_{\alpha} , a+\alpha x_{\alpha}]}(y) - y\, 1_{[-a+\alpha_0 x_{\alpha},-a+\alpha x_{\alpha}]}(y) .
$$
From (\ref{noyau_AR}), we have 
\begin{subequations}
\begin{eqnarray*}
& & (P_\alpha f_\alpha)(x_\alpha) = \int_{0}^{a} (y+\alpha x_\alpha)\, \nu(y)dy - \int_{-2a}^{-a} (y+\alpha x_\alpha)\, \nu(y)dy \\
& & (P_{\alpha_0} f_\alpha)(x_\alpha) = \int_{a}^{2a} (y+\alpha_0 x_\alpha)\, \nu(y)dy - \int_{-a}^{0} (y+\alpha_0 x_\alpha)\, \nu(y)dy.
\end{eqnarray*}
\end{subequations}
Set $J:= \int_{-a}^a y\, \nu(y)dy - \int_{-2a}^{-a} y\, \nu(y)dy - \int_a^{2a}  y\, \nu(y)dy$. We obtain 
$$(P_\alpha f_\alpha)(x_\alpha) - (P_{\alpha_0} f_\alpha)(x_\alpha) =  J + \alpha x_\alpha\bigg(\int_0^a \nu(y)dy - \int_{-2a}^{-a} \nu(y)dy\bigg) + \alpha_0 x_\alpha\bigg(\int_{-a}^0 \nu(y)dy - \int_a^{2a} \nu(y)dy\bigg).$$
Since $V(x_\alpha) = 1 + x_\alpha$ and $x_\alpha\r+\infty$ as $\alpha\downarrow\alpha_0$, it follows that  
$$\lim_{\alpha\downarrow\alpha_0} \frac{(P_\alpha f_\alpha)(x_\alpha) - (P_{\alpha_0} f_\alpha)(x_\alpha)}{V(x_\alpha)} = \alpha_0 \bigg(\int_{-a}^a \nu(y)dy - \int_{-2a}^{-a} \nu(y)dy  - \int_a^{2a} \nu(y)dy\bigg) \neq 0.$$
Finally, from $$\|P_\alpha - P_{\alpha_0}\|_{1}= \|P_\alpha - P_{\alpha_0}\|_{{\cal B}_V} = \sup_{|f|\leq V}\sup_{x\in\R}\frac{\big|(P_\alpha f)(x) - (P_{\alpha_0} f)(x)\big|}{V(x)} \geq 
\frac{\big|(P_\alpha f_\alpha)(x_\alpha) - (P_{\alpha_0} f_\alpha)(x_\alpha)\big|}{V(x_\alpha)},$$
we deduce that $\|P_\alpha - P_{\alpha_0}\|_{{\cal B}_V}$ does not go to 0 when $\alpha\downarrow\alpha_0$. 
\end{ex}

Reinforcing the continuity condition~(\ref{Cont_01}) of Theorem~\ref{pro-bv} allows us to obtain the following refinement. 
\begin{cor} \label{drift}
Assume that Conditions~{\em (\ref{V-geo-cond})} and {\em (\ref{drift-gene-cond})} hold and that there exists $\beta\in(0,1]$ such that\footnote{Under the drift  condition,  (\ref{V-geo-cond}) implies (\ref{V-geo-cond-beta}) by Jensen's inequality for $\beta\in (0,1]$.} 
\begin{equation} \label{V-geo-cond-beta} 
\exists \kappa_\beta\in(0,1),\quad \|P_0^n - \pi_0(\cdot)1_\X\|_{\beta} = O(\kappa_\beta^n) \tag{$V_\beta$}
\end{equation}
\begin{equation} \label{Cont_beta_1}
\text{and }\lim_{\varepsilon\r0} \|P_\varepsilon - P_0\|_{\beta,1} = 0. 
\end{equation}
Then there exist $\varepsilon_1\in(0,\varepsilon_0]$ (from Theorem~\ref{pro-bv}) and $C\in(0,+\infty)$ such that: 
$$\forall \varepsilon\in(-\varepsilon_1,\varepsilon_1),\quad \sup_{\|f\|_\beta\leq1}\big|\pi_\varepsilon(f) - \pi_0(f)\big| \leq C\, \|P_\varepsilon - P_0\|_{\beta,1}.$$
\end{cor}

\begin{rem} \label{Rem_papier_qc}
Let $P$ be any Markov kernel on $\X$. If we have $P^N V \leq \delta^N V + L\, 1_{\X}$ for some $N\in\N^*$, $L\in(0,+\infty)$ and $\delta\in(0,1)$ and if $P^\ell$ is compact from $\cB_0$ to $\cB_V$ (for some $\ell\geq1$), then $P$ is a power-bounded quasi-compact operator on $\cB_V$. Moreover its essential spectral radius is such that  $r_{ess}(P)\leq\delta$.  Actually, denoting by $\delta_V(P)$ the infinum of the real numbers $\delta$ for which the above inequality holds for some $N\in\N^*$ and $L\in(0,+\infty)$, we have $r_{ess}(P)\leq\delta_V(P)$, with equality in many cases. This follows from the Doeblin-Fortet inequalities (\ref{doeb-fortet-dual}) and \cite{hen}, since $(P^*)^\ell$ is compact from $\cB_V'$ to $\cB_0'$. If in addition $P$ satisfies  aperiodicity and irreducibility conditions, then $P$ is $V$-geometrically ergodic. Such results are fully detailed and applied in \cite{djl-qc} (without any perturbation issues). Now, if $\rho_V(P)$ denotes the infinum bound of the real numbers $\kappa_1$ such that {\em (\ref{V-geo-cond})} holds, then we have $\rho_V(P)\geq r_{ess}(P)$. If there are eigenvalues $\lambda$ of $P$ on $\cB_V$ such that $r_{ess}(P)<|\lambda|<1$, then $\rho_V(P)$ is the maximal modulus of such eigenvalues (which are in finite number from the definition of $r_{ess}(P)$); if not, we have $\rho_V(P) = r_{ess}(P)$. 
For instance, if $(X_n)_{n\in\N}$ is defined by (\ref{def-AR}) and if the probability density function  $\nu(\cdot)$ of the noise has a moment of order $r\in [1,+\infty)$, then $(X_n)_{n\in\N}$ is $(1+|\cdot|)^r$-geometrically ergodic with $r_{ess}(P)\leq |\alpha|^r$ and $\rho_V(P)=|\alpha|$. See \cite[Sec. 8]{Wu04} or \cite[Sec.~5]{djl-qc} for details. 
\end{rem}

Now we are interested in asymptotic expansions of the perturbed invariant distributions. Under Condition~(\ref{VG}) and the continuity condition~(\ref{Cont_11}), expansions related to the generalized potential $R = (I-P_0+\Pi_0)^{-1}$,  with $\Pi_0(\cdot) = \pi_0(\cdot)1_\X$, are given in \cite{kar81,kar96,hei-hor}, namely
$$\pi_\varepsilon = \pi_0(I-D_\varepsilon R)^{-1} = \sum_{k\geq0} \pi_0(D_\varepsilon R)^k = \pi_0 + \sum_{k=1}^n \pi_0(D_\varepsilon R)^k + O((D_\varepsilon R)^{n+1}),$$
when $D_\varepsilon = P_\varepsilon  - P_0$ is such that $\|D_\varepsilon R\|_1 < 1$. In general this is not a Taylor expansion, except for instance when $D_\varepsilon$ has the form $D_\varepsilon = \epsilon D$ with $D\in\cL(\cB_1)$. This special case is discussed in \cite{AAN04}. In fact, for general perturbations, even in the case when (\ref{Cont_11}) is fulfilled, obtaining Taylor expansions for $\pi_\varepsilon$ causes difficulties when the derivatives of the perturbed kernels $P_\varepsilon(x,\cdot)$ w.r.t.~to $\varepsilon$ yield some weights (see the term $x^k$ in (\ref{deri-formel}) for order of regularity $k\geq1$). For this question, the derivation procedure (see \cite{gouliv,fl}) based on the Keller-Liverani theorem is of great interest. 

To illustrate this approach, we only consider the special instance of autoregressive model  introduced in Example~\ref{AR_weak_pasStrong}. More specifically, we assume that $(X_n)_{n\in\N}$ is given by (\ref{def-AR}) and that the probability density function $\nu(\cdot)$ of the noise has a moment of order $r$, $\int |x|^r \nu(x) dx <\infty$, for some  $r\in [1,+\infty)$. Set $V(x) := (1+|x|)^r$, $x\in\R$. Then, for each $\alpha\in(-1,1)$, $(X_n)_{n\in\N}$ is $V$-geometrically ergodic with an invariant distribution $\pi_\alpha$. Moreover we assume that $r$ is not an integer (ie.~$r>{\lfloor r\rfloor}$ where $\lfloor \cdot \rfloor$ stands for the integer part function on $\R$), that $\nu(\cdot)$ is positive and ${\lfloor r\rfloor}+1$-times continuously differentiable on $\R$, with 
\begin{equation} \label{cond-derive-nu} 
j=1,\ldots,{\lfloor r\rfloor}+1, \quad \sup_{t\in\R}\frac{|\nu^{(j)}(t)|}{\nu(t)} < \infty.
\end{equation}
Finally suppose that for all $x_0\in \R$, there exist a neighborhood $V_{x_0}$ of $x_0$ and a non-negative measurable function $q_{x_0}(\cdot)$ such that $\int_\R (1+|y|)^{r}\, q_{x_0}(y)\, dy < \infty$ and: $\forall y\in\R,\ \forall v\in V_{x_0},\ \nu(y+v) \leq q_{x_0}(y)$. 
\begin{pro} \label{auto-reg} 
Under the previous assumptions, there exists $\beta\equiv\beta(r)\in(0,1)$ such that the map $\alpha\mapsto \pi_\alpha$ is ${\lfloor r\rfloor}$-times continuously differentiable from $(-1,1)$ to the dual space $\cB_{\beta}'$ of $\cB_{\beta}$. In particular, for all $\alpha\in(-1,1)$, there exist $\lfloor r\rfloor$ signed measures $\mu_{\alpha,1},\ldots,\mu_{\alpha,{\lfloor r\rfloor}}$ on $\R$ such that: 
$$\forall A\in\cX,\quad \pi_{\alpha+\varepsilon}(A) = \pi_{\alpha}(A) + \sum_{j=1}^{\lfloor r\rfloor} \frac{\varepsilon^j}{j!}\mu_{\alpha,j}(A) + \varepsilon^{\lfloor r\rfloor}R_{\varepsilon}(A) \quad \text{with } \lim_{\varepsilon\r0}\sup_{A\in{\cal X}}\big|R_{\varepsilon}(A)\big| = 0.$$
\end{pro}

The previous conditions (\ref{cond-derive-nu}) are well suited to densities of the form $\nu(x) = c(1+|x|)^{-\gamma}$. For other forms of densities, they can be easily adapted in order that Proposition~\ref{auto-reg} works. 

\section{Proofs}

\begin{proof}{ of Theorem~\ref{pro-bv}} 
 We shall repeatedly use the fact that a bounded linear operator (between two normed vector spaces) and its adjoint have the same norm. The adjoint of the rank-one projection $\Pi_0(f) = \pi_0(f)1_\X$ is defined by: $\forall f'\in\cB_1',\ \Pi_0^*(f') :=  f'(1_\X)\, \pi_0$. 

First, Assumption~(\ref{V-geo-cond}) gives 
$$\|(P_0^*)^n - \Pi_0^*\|_{{\cal B}_1'} = O(\kappa_1^n).$$
Second, it follows from the continuity assumption (\ref{Cont_01}) of Theorem~\ref{pro-bv} that
 $$\lim_{\varepsilon\r0} \|P_\varepsilon^*  - P_0^*\|_{{\cal B}_1',{\cal B}_0'} = 0,$$ 
Third, Assumption~(\ref{drift-gene-cond}) gives: 
\begin{equation} \label{doeb-fortet-dual}
\forall \varepsilon\in(-\varepsilon_0,\varepsilon_0),\ \forall f'\in\cB_1',\quad \|P_\varepsilon^{*N}f'\|_{{\cal B}'_1} \leq \delta^N\|f'\|_{{\cal B}'_1} + L \|f'\|_{{\cal B}'_0}.
\end{equation} 
Indeed, recall that $\cB_1'$ and $\cB_0'$ are Banach lattices and that, for each $g'\in\cB_1'$, we have $\|g'\|_{{\cal B}_0'} = \|\, |g'|\, \|_{{\cal B}_0'} = \langle |g'|,1_{\X}\rangle$ and $\|g'\|_{{\cal B}_1'} =  \|\, |g'|\, \|_{{\cal B}_1'} = \langle |g'|,V\rangle$, where $\langle \cdot,\cdot\rangle$ stands for the duality brackets in both $\cB_0'\times \cB_0$ and $\cB_1'\times \cB_1$. Next, observe that, for every $f'\in\cB_1'$ and every $f\in\cB_1$ such that $|f|\leq V$, we have $|\langle (P^*)^Nf',f\rangle| \leq \langle |f'|, P^NV\rangle$. Hence (\ref{drift-gene-cond}) gives (\ref{doeb-fortet-dual}) as claimed.

The three previous facts show that $\{P_\varepsilon^*\}_{|\varepsilon|<\varepsilon_0}$ satisfies  the assumptions of \cite{liv04} on $\cB_1'$ (see also \cite{fl} for the use of the perturbation Keller-Liverani theorem in a Markov context). Therefore, for all $\kappa\in(\widehat\kappa,1)$, there exists $\varepsilon_1\in(0,\varepsilon_0]$ such that, for all  $\varepsilon\in(-\varepsilon_1,\varepsilon_1)$, the following properties hold: 
there exist $\lambda_\varepsilon\in\C$ satisfying $\lim_{\varepsilon\r0} \lambda_\varepsilon = 1$ and a rank-one projection $\Pi_\varepsilon'$ on $\cB_1'$ such that 
\begin{equation} \label{conseq-K-L}
\sup_{|\varepsilon|<\varepsilon_1}\|(P_\varepsilon^*)^n - \lambda_\varepsilon^n\Pi_\varepsilon'\|_{{\cal B}_1'} = O(\kappa^n) \quad \text{ and }\quad \lim_{\varepsilon\r0} \|\Pi_\varepsilon' - \Pi_0^*\|_{{\cal B}_1',{\cal B}_0'} = 0. 
\end{equation} 
Up to reduce $\varepsilon_1$, we have $\lambda_\varepsilon = 1$. Indeed one may assume that $|\lambda_\varepsilon|>\kappa$. Then we deduce from $1_{\X}=P_\varepsilon^n1_{\X}$ and (\ref{conseq-K-L}) that, for any $f'\in\cB_1'$, we have $\lim_n \lambda_\varepsilon^{-n}\langle f', 1_\X\rangle = 
\lim_n \lambda_\varepsilon^{-n} \langle (P_\varepsilon^*)^n f', 1_\X \rangle = \langle \Pi_\varepsilon'f', 1_\X\rangle$. Thus $(\lambda_\varepsilon^{-n})_{n\in\N}$ converges in $\C$, so that we have, either $\lambda_\varepsilon=1$, or $|\lambda_\varepsilon| >1$. Moreover the sequence $(P_\varepsilon^n)_{n\in\N}$ is bounded in $\cL(\cB_V)$ from (\ref{drift-gene-cond}). Thus $((P_\varepsilon^*)^n)_{n\in\N}$ is bounded in $\cL(\cB_V')$, and (\ref{conseq-K-L}) then implies that $(\lambda_\varepsilon^{n})_{n\in\N}$ is bounded in $\C$. Thus $\lambda_\varepsilon = 1$. 

Hence $((P_\varepsilon^*)^n)_{n\in\N}$ is Cauchy, and so is $(P_\varepsilon^n)_{n\in\N}$ in $\cL(\cB_1)$. Therefore 
\begin{equation} \label{Cv_proj}
\exists\, \Pi_\varepsilon\in\cL(\cB_1),\quad \lim_n P_\varepsilon^n = \Pi_\varepsilon \text{ in } \cL(\cB_1).
\end{equation}
 Thus $(P_\varepsilon^n)^*\r \Pi_\varepsilon^*$ in $\cL(\cB_1')$, and so $\Pi_\varepsilon' = \Pi_\varepsilon^*$. We have obtained that 
 $$\sup_{|\varepsilon|<\varepsilon_1}\|P_\varepsilon^n - \Pi_\varepsilon\|_{1} = O(\kappa^n) \quad \text{ and } \quad
\lim_{\varepsilon\r0}\|\Pi_\varepsilon - \Pi_0\|_{0,1}=0.$$ 
From (\ref{Cv_proj}) it follows that $\Pi_\varepsilon$ is a positive projection on $\cB_1$ satisfying $\Pi_\varepsilon P_\varepsilon = P_\varepsilon \Pi_\varepsilon = \Pi_\varepsilon$. 
Let us prove that $\Pi_\varepsilon$ is rank-one. We know that $\Pi_\varepsilon'$ is rank-one, namely $\Pi_\varepsilon'(\cdot) := \langle \phi_\varepsilon,\cdot \rangle \, a_\varepsilon'$ for some $\phi_\varepsilon\neq0$ in the dual space of $\cB_1'$ and for some $a_\varepsilon'\neq0$ in $\cB_1'$. From  $\Pi_\varepsilon^* = \Pi_\varepsilon'$, we obtain $\ker \Pi_\varepsilon = \ker a_\varepsilon'$. Thus $\ker\Pi_\varepsilon$ is of codimension one, so that we have $\dim\Im\Pi_\varepsilon=1$ (use $\cB_1 = \ker\Pi_\varepsilon \oplus \Im\Pi_\varepsilon$). 

The last fact shows that $\Pi_\varepsilon(\cdot) := e_\varepsilon'(\cdot)\, 1_\X$ for some $e_\varepsilon'\geq0$ in $\cB_1'$. Since $\forall A\in\cX,\ \lim_n P_\varepsilon^n(x_0,A) = e_\varepsilon'(1_A)$ (for some  fixed $x_0\in\X$) and $P_\varepsilon^n(x_0,\cdot)$ is a probability measure on $(\X,\cX)$ for each $n\geq1$, we deduce from the Vitali-Hahn-Saks theorem that $A\mapsto \pi_\varepsilon(A) := e'_\varepsilon(1_A)$ is a probability measure on $(\X,\cX)$. Clearly $\pi_\varepsilon$ is $P_\varepsilon$-invariant. Next $\pi_\varepsilon$ and $e'_\varepsilon$ coincide on $\cB_0$: indeed every $f\in\cB_0$ can be approached uniformly on $\X$ by a sequence $(f_n)_n$ of simple functions, so that $\pi_\varepsilon(f) = \lim_n\pi_\varepsilon(f_n) = \lim_ne'_\varepsilon(f_n) = e'_\varepsilon(f)$ (the last convergence holds since $e'_\varepsilon\in\cB_1'$). We obtain $\pi_\varepsilon(V1_{[V\leq n]}) = e'_\varepsilon(V1_{[V\leq n]}) \leq e'_\varepsilon(V) < \infty$ for every $n\geq1$ since $e_\varepsilon'\geq0$. 
The monotone convergence theorem gives $\pi_\varepsilon(V) < \infty$. Thus $f\mapsto\pi_\varepsilon(f)$ is in $\cB_1'$. Since $\forall f\in\cB_1,\ \lim_n P_\varepsilon^nf  = e_\varepsilon'(f)\, 1_\X$ in $\cB_1$, 
we deduce that $e_\varepsilon'=\pi_\varepsilon$ on $\cB_1$ from the invariance of $\pi_\varepsilon$. Finally the last assertion of Theorem~\ref{pro-bv} follows from  $\lim_{\varepsilon\r0}\|\Pi_\varepsilon - \Pi_0\|_{0,1}=0$. 
 \end{proof}

\noindent\begin{proof}{ of Corollary~\ref{drift}} 
Let $\kappa\in(\max(\kappa_\beta,\kappa_1),1)$, and let $\Gamma$ be the oriented circle in $\C$, centered at $z=1$ and with radius less than $(1-\kappa)/2$. Note that the assumptions of Corollary~\ref{drift} imply those of Theorem~\ref{pro-bv}. We know from  \cite{liv04} and the proof of Theorem~\ref{pro-bv} that 
$M := \sup_{z\in\Gamma,|\varepsilon|<\varepsilon_1}\|(zI- P_\varepsilon)^{-1}\|_{1} < \infty$.  
 Moreover, the rank-one eigenprojection $\Pi_\varepsilon$ is from the standard spectral theory: 
\begin{equation} \label{proj-int}
\Pi_\varepsilon = \frac{1}{2i\pi}\oint_{\Gamma} (zI- P_\varepsilon)^{-1} dz.
\end{equation}
Next (\ref{V-geo-cond-beta}) gives  $M_0 := \sup_{z\in\Gamma}\|(zI- P_0)^{-1}\|_{\beta} < \infty$.  Assume that $\|1_\X\|_1=1$ (to simplify). Then, from 
 $$(zI- P_\varepsilon)^{-1}  - (zI- P_0)^{-1} = (zI- P_\varepsilon)^{-1} (P_\varepsilon - P_0) (zI- P_0)^{-1},$$
  we obtain  that for every  $f\in\cB_\beta$
  $$\big|\pi_\varepsilon(f) - \pi_0(f)\big| = \|\Pi_\varepsilon(f) - \Pi_0(f)\|_1 
\leq \frac{(1-\kappa)}{2}\, M\, \|P_\varepsilon - P_0\|_{\beta,1}\, M_0\, \|f\|_\beta.$$ 
\end{proof}
\begin{proof}{ of Proposition~\ref{auto-reg}} 
Recall that $V(x) := (1+|x|)^r$, $x\in\R$. Let $a_0\in(0,1)$.
We apply the derivation procedure of \cite[App.~A]{fl} to the family $\{P_\alpha,\, \alpha\in(-a_0,a_0)\}$ with respect to the Banach spaces $\cB_\beta$. Let $\beta\in(0,1]$. For any $\alpha\in(-a_0,a_0)$, the Markov kernel $P_\alpha$ of the autoregressive model $(X_n)_{n\in\N}$ satisfies the drift condition and the aperiodicity/irreducibility assumptions of \cite{mey}. Hence each $P_\alpha$ satisfies~(\ref{V-geo-cond-beta}). In the same way, $\{P_\alpha\}_{|\alpha|<a_0}$ satisfies (\ref{drift-gene-cond}), so that it also satisfies (\ref{drift-gene-cond}) with respect to the function $V^\beta$ thanks to Jensen's inequality. Let $\alpha\in(-a_0,a_0)$. Then 
$\lim_{\varepsilon\r0} \|P_{\alpha+\varepsilon} - P_{\alpha}\|_{0,\beta} = 0$ from Lemma~\ref{lem-cont-AR} below. Hence the perturbation Keller-Liverani theorem can be applied in $\cL(\cB_\beta)$ to the family $(P_{\alpha+\epsilon})_\epsilon$ (as seen with $\beta=1$ in the proof of Theorem~\ref{pro-bv}). This gives the following spectral properties: there exists $\kappa_{\alpha,\beta}\in(0,1)$ such that, for all $\kappa\in(\kappa_{\alpha,\beta},1)$, there exists $\widehat\varepsilon>0$ such that the resolvents $(zI-P_{\alpha+\varepsilon})^{-1}$ are well-defined and uniformly bounded in $\cL(\cB_\beta)$ provided that $|\varepsilon|\leq\widehat\varepsilon$ and $z\in\C$ satisfies  $|z|\geq\kappa$ and $|z-1|\geq(1-\kappa)/2$. 

Now let us introduce the formal derivative operators of $\alpha\mapsto P_\alpha$: 
\begin{equation} \label{deri-formel}
\forall x\in\R,\quad \big(P_{k,\alpha} f\big)(x) = (-1)^k x^k \int_\R f(y) \nu^{(k)}(y-\alpha x)\, dy. 
\end{equation}
\begin{lem} \label{lem-cont-AR}
Let $(\beta,\beta')\in[0,1]^2$. 
\begin{enumerate}[(i)]
	\item If $\beta+k/r<\beta'\leq1$ for $k\in\{0,\ldots,\lfloor r\rfloor\}$, then $\alpha \mapsto P_{k,\alpha}$ is continuous from $(-a_0,a_0)$ to $\cL(\cB_\beta,\cB_{\beta'})$. 
	\item If $\beta+(k+1)/r< \beta'\leq1$ for $k\in\{0,\ldots,\lfloor r\rfloor-1\}$, then $\alpha \mapsto P_{k,\alpha}$ is continuously differentiable from $(-a_0,a_0)$ to $\cL(\cB_\beta,\cB_{\beta'})$, with $\frac{dP_{k,\alpha}}{d\alpha}(\alpha) = P_{k+1,\alpha}$. 
	\end{enumerate}
\end{lem}

In a first stage, let us admit that Lemma~\ref{lem-cont-AR} holds and apply \cite[App.~A]{fl}. For $\beta\in(0,1)$ and $\sigma>0$, we set: $T_0(\beta) = \beta + \sigma/r$ and $T_1(\beta) = \beta + (1+\sigma)/r$. Lemma~\ref{lem-cont-AR} gives for $j=1,\ldots,\lfloor r\rfloor$:  
\begin{itemize}
\item[--] if $T_0(\beta) \in(0,1]$, then the map $\alpha \mapsto P_\alpha$ is continuous from $(-a_0,a_0)$ to $\cL(\cB_\beta,\cB_{T_0(\beta)})$,
\item[--] if $T_1(T_0 T_1)^{j-1}(\beta) \in (0,1]$, then the map $\alpha \mapsto P_\alpha$ is $j$-times continuously differentiable from  $(-a_0,a_0)$ to $\cL(\cB_\beta,\cB_{T_1(T_0 T_1)^{j-1}(\beta)})$. 
\end{itemize}
Let $\beta_r \in (0,1 - \lfloor r\rfloor/r)$, and let $\sigma >0$ such that $\beta_r + [(2\lfloor r\rfloor+1) \sigma  + \lfloor r\rfloor]/r = 1$. In other words we have $(T_0T_1)^{\lfloor r\rfloor} T_0(\beta_r) = 1$. Let $\alpha_0\in(-a_0,a_0)$. Then it follows from \cite[Prop.~A.1]{fl} that there exists $\tilde\kappa\in(0,1)$ such that, for all  $z\in\cD := \{z\in\C : |z|\geq\tilde\kappa,\, |z-1|\geq(1-\tilde\kappa)/2\}$, the map $\alpha\mapsto(zI-P_\alpha)^{-1}$ is ${\lfloor r\rfloor}$-times continuously differentiable from some open interval $I_{\alpha_0}$  centered at $\alpha_0$ into $\cL(\cB_{\beta_r},\cB_1)$. 
Moreover the derivatives (up to the order ${\lfloor r\rfloor}$) of the last map are uniformly bounded in $(\alpha,z)\in I_{\alpha_0}\times\cD$. Next, let us define $\Pi_\alpha(f) = \pi_\alpha(\cdot)1_\R$ for any $f\in\cB_{\beta_r}$. From the previous facts and standard spectral calculus (see (\ref{proj-int})), the map $\alpha \mapsto \Pi_\alpha(f)$ is ${\lfloor r\rfloor}$-times continuously differentiable from $I_{\alpha_0}$ to $\cL(\cB_{\beta_r},\cB_1)$, so is $\alpha \mapsto \pi_\alpha(\cdot)$ from $I_{\alpha_0}$ to $\cB_{\beta_r}'$. Since $\alpha_0$ is any element in $(-a_0,a_0)$, with arbitrary $a_0\in(-1,1)$, this proves the first assertion of Proposition~\ref{auto-reg}. 

Finally let $\alpha\in(-1,1)$. For all borel set $A$ of $\R$, we have $\frac{d\pi_\alpha}{d\alpha}(A) = \lim_{h\r0}(\pi_{\alpha+h}(A) - \pi_\alpha(A))/h$. Since $\pi_{\alpha+h}(\cdot) - \pi_\alpha(\cdot)$ is a (signed) measure on $\R$, it follows from the Vitali-Hahn-Saks theorem that there exists a (signed) measure $\mu_{\alpha,1}$ on $\R$ such that we have: $\frac{d\pi_\alpha}{d\alpha}(A) = \mu_{\alpha,1}(A)$. An obvious induction gives the same conclusion for the derivatives of order $j=2,\ldots,\lfloor r\rfloor$. 
\end{proof}
\begin{proof}{ of Lemmas~\ref{lem-cont-AR}} 
Set $A_k := \sup_{t\in\R} |\nu^{(k)}(t)|/\nu(t)$ and $A:=\max_{0\leq k \leq \lfloor r\rfloor +1}\, A_k$. 
Note that there exists $B\equiv B(\beta)$ such that, for all $\alpha\in[-a_0,a_0]$, we have $P_\alpha V^\beta \leq B V^\beta$. 

For $\alpha\in(-a_0,a_0)$, $\beta\in(0,1]$, $f\in\cB_\beta$ and $x\in\R$, we denote $F_k(\alpha) := (P_{k,\alpha} f)(x)$. 
From the assumptions on $\nu$ and Lebesgue's theorem, for every $k=0,\ldots,\lfloor r\rfloor$, $F_k$ is differentiable on $[-a_0,a_0]$, with 
$\frac{\partial F_k}{\partial \alpha}(\alpha) = F_{k+1}(\alpha)$. 

Let $k\in\{0,\ldots,\lfloor r\rfloor\}$, $(\beta,\beta')\in[0,1]^2$ such that $\beta+k/r<\beta'\leq1$, and let  $0<\sigma \leq 1$ be such that $\beta + (k+\sigma)/r = \beta'$. For any $(\alpha,\alpha')\in(-a_0,a_0)^2$, we obtain (use Taylor expansion of $F_k$ for (\ref{eq2})): 
\begin{subequations}
\begin{eqnarray}
 \big|F_k(\alpha) - F_k(\alpha')\big| & \leq & 2\, A_{k} B\, \|f\|_\beta\, V(x)^{\beta+\frac{k}{r}}.
  \label{eq1} \\
\big|F_k(\alpha) - F_k(\alpha')\big| & \leq & |\alpha-\alpha'|\, A_{k+1} B\, \|f\|_\beta\, V(x)^{\beta+\frac{k+1}{r}}.
\label{eq2}
\end{eqnarray}
\end{subequations}
Multiplying (\ref{eq1}) (to the power $1-\sigma$) by (\ref{eq2}) (to the power $\sigma$) gives: 
$$|F_k(\alpha) - F_k(\alpha')| \leq 2|\alpha-\alpha'|^\sigma\, A B\, \|f\|_\beta\, V(x)^{\beta+\frac{k+\sigma}{r}},$$ 
which rewrites as: $\|P_{k,\alpha} f - P_{k,\alpha'}f\|_{\beta'} \leq 2AB\, |\alpha-\alpha'|^\sigma \, \|f\|_\beta$. 
The first assertion of Lemma~\ref{lem-cont-AR} is proved. 

Next let $k\in\{0,\ldots,\lfloor r\rfloor-1\}$, $(\beta,\beta')\in[0,1]^2$ such that $\beta+(k+1)/r< \beta'\leq1$, and let $0<\sigma \leq 1$ be such that $\beta + (k+1+\sigma)/r = \beta'$. For any $(\alpha,\alpha')\in(-a_0,a_0)^2$, we obtain:  
\begin{subequations}
\begin{eqnarray}
& &  \big|F_k(\alpha) - F_k(\alpha')-(\alpha-\alpha')F_{k+1}(\alpha')\big| \leq 2\,|\alpha-\alpha'| \, A_{k+1} B\, \|f\|_\beta\, V(x)^{\beta+\frac{k+1}{r}}
\label{eq3} \\
& & \big|F_k(\alpha) - F_k(\alpha')-(\alpha-\alpha')F_{k+1}(\alpha')\big| \leq \frac{|\alpha-\alpha'|^2}{2}\, A_{k+2} B\, \|f\|_\beta\, V(x)^{\beta+\frac{k+2}{r}}. \label{eq4}
\end{eqnarray}
\end{subequations}
Then, by multiplying (\ref{eq3}) (with the power $1-\sigma$) and (\ref{eq4}) (with the power $\sigma$), we obtain: 
$$|F_k(\alpha) - F_k(\alpha')-(\alpha-\alpha')F_{k+1}(\alpha')| \leq 2\, |\alpha-\alpha'|^{1+\sigma}\, A B\, \|f\|_\beta\, V(x)^{\beta+\frac{k+1+\sigma}{r}},$$ 
which rewrites as: $\|P_{k,\alpha} f - P_{k,\alpha'}f-(\alpha-\alpha') P_{k+1,\alpha'}f\|_{\beta'} \leq 2 A B\,  |\alpha-\alpha'|^{1+\sigma} \, \|f\|_\beta$. 
\end{proof}


\end{document}